\documentclass[11pt]{amsart}
\usepackage{amssymb,mathrsfs}
\title{Invariant Metrics and Laplacians on Siegel-Jacobi Disk}

\author{Jae-Hyun Yang}
\address{Department of Mathematics, Inha University,
Incheon 402-751, Republic of Korea}
\email{jhyang@inha.ac.kr }

\begin{document}

\thanks{
\indent This work was supported by Inha University Research Grant}

\begin{abstract} Let ${\mathbb D}_n $ be the generalized unit disk of degree $n$.
In this paper, we find Riemannian metrics on the Siegel-Jacobi
disk ${\mathbb D}_n \times {\mathbb C}^{(m,n)}$
 which are invariant under the natural action of the Jacobi group
explicitly and also compute the Laplacians of these invariant
metrics explicitly. These are expressed in terms of the trace
form. We give a brief remark on the theory of harmonic analysis on
the Siegel-Jacobi disk ${\mathbb D}_n \times {\mathbb C}^{(m,n)}$.
\end{abstract}

\maketitle

\newtheorem{theorem}{Theorem}[section]
\newtheorem{lemma}{Lemma}[section]
\newtheorem{proposition}{Proposition}[section]
\newtheorem{remark}{Remark}[section]
\newtheorem{definition}{Definition}[section]

\renewcommand{\theequation}{\thesection.\arabic{equation}}
\renewcommand{\thetheorem}{\thesection.\arabic{theorem}}
\renewcommand{\thelemma}{\thesection.\arabic{lemma}}
\newcommand{\bbr}{\mathbb R}
\newcommand{\bbs}{\mathbb S}
\newcommand{\bn}{\bf n}
\def\charf {\mbox{{\text 1}\kern-.24em {\text l}}}

\newcommand\BD{\mathbb D}
\newcommand\BH{\mathbb H}
\newcommand\BR{\mathbb R}
\newcommand\BC{\mathbb C}
\newcommand\lrt{\longrightarrow}
\newcommand\lmt{\longmapsto}
\newcommand\CX{{\Cal X}}
\newcommand\td{\bigtriangledown}
\newcommand\pdx{ {{\partial}\over{\partial x}} }
\newcommand\pdy{ {{\partial}\over{\partial y}} }
\newcommand\pdu{ {{\partial}\over{\partial u}} }
\newcommand\pdv{ {{\partial}\over{\partial v}} }
\newcommand\PZ{ {{\partial}\over {\partial Z}} }
\newcommand\PW{ {{\partial}\over {\partial W}} }
\newcommand\PZB{ {{\partial}\over {\partial{\overline Z}}} }
\newcommand\PWB{ {{\partial}\over {\partial{\overline W}}} }
\newcommand\PX{ {{\partial}\over{\partial X}} }
\newcommand\PY{ {{\partial}\over {\partial Y}} }
\newcommand\PU{ {{\partial}\over{\partial U}} }
\newcommand\PV{ {{\partial}\over{\partial V}} }
\renewcommand\th{\theta}
\renewcommand\l{\lambda}
\renewcommand\k{\kappa}
\newcommand\G{\Gamma}
\newcommand\Om{\Omega}
\newcommand\s{\sigma}
\newcommand\g{\gamma}
\newcommand\PWE{ \frac{\partial}{\partial \overline W}}

\newcommand\PO{ {{\partial}\over {\partial \Omega}} }
\newcommand\PE{ {{\partial}\over {\partial \eta}} }
\newcommand\POB{ {{\partial}\over {\partial{\overline \Omega}}} }
\newcommand\PEB{ {{\partial}\over {\partial{\overline \eta}}} }

%$^t\!\left(\PWE\right)$

\begin{section}{{\bf Introduction}}
\setcounter{equation}{0}
\end{section}

For a given fixed positive integer $n$, we let
$${\mathbb H}_n=\,\{\,\Omega\in \BC^{(n,n)}\,|\ \Om=\,^t\Om,\ \ \ \text{Im}\,\Om>0\,\}$$
be the Siegel upper half plane of degree $n$ and let
$$Sp(n,\BR)=\{ M\in \BR^{(2n,2n)}\ \vert \ ^t\!MJ_nM= J_n\ \}$$
be the symplectic group of degree $n$, where
$$J_n=\begin{pmatrix} 0&I_n \\
                   -I_n&0 \\ \end{pmatrix}.$$
We see that $Sp(n,\BR)$ acts on $\BH_n$ transitively by
\begin{equation}M\cdot\Om=(A\Om+B)(C\Om+D)^{-1}, \end{equation}
where $M=\begin{pmatrix} A&B\\ C&D\end{pmatrix}\in Sp(n,\BR)$ and
$\Om\in \BH_n.$

For two positive integers $n$ and $m$, we consider the Heisenberg
group
$$H_{\BR}^{(n,m)}=\{\,(\l,\mu;\k)\,|\ \l,\mu\in \BR^{(m,n)},\ \k\in \BR^{(m,m)},\ \
\k+\mu\,^t\l\ \text{symmetric}\ \}$$ endowed with the following
multiplication law
$$(\l,\mu;\k)\circ (\l',\mu';\k')=(\l+\l',\mu+\mu';\k+\k'+\l\,^t\mu'-
\mu\,^t\l').$$ We define the semidirect product of $Sp(n,\BR)$ and
$H_{\BR}^{(n,m)}$
$$G^J:=Sp(n,\BR)\ltimes H_{\BR}^{(n,m)}$$
endowed with the following multiplication law
$$\Big(M,(\lambda,\mu;\kappa)\Big)\cdot\Big(M',(\lambda',\mu';\kappa')\Big) =\,
\Big(MM',(\tilde{\lambda}+\lambda',\tilde{\mu}+ \mu';
\kappa+\kappa'+\tilde{\lambda}\,^t\!\mu'
-\tilde{\mu}\,^t\!\lambda')\Big)$$ with $M,M'\in Sp(n,\BR),
(\lambda,\mu;\kappa),\,(\lambda',\mu';\kappa') \in
H_{\BR}^{(n,m)}$ and
$(\tilde{\lambda},\tilde{\mu})=(\lambda,\mu)M'$. We call this
group $G^J$ the Jacobi group of degree $n$ and index $m$. Then we
get the {\it natural action} of $G^J$ on $\BH_n\times
\BC^{(m,n)}$\,(cf.\,[1-2],\,[7-9],\,[14]) defined by
\begin{equation}\Big(M,(\lambda,\mu;\kappa)\Big)\cdot (\Om,Z)=\Big(M\cdot\Om,(Z+\lambda \Om+\mu)
(C\Om+D)^{-1}\Big), \end{equation}
where $M=\begin{pmatrix} A&B\\
C&D\end{pmatrix} \in Sp(n,\BR),\ (\lambda,\mu; \kappa)\in
H_{\BR}^{(n,m)}$ and $(\Om,Z)\in \BH_n\times \BC^{(m,n)}.$ We note
that the action (1.2) is transitive.

\vskip 0.31cm
\newcommand\bw{d{\overline W}}
\newcommand\bz{d{\overline Z}}
\newcommand\be{d{\overline \eta}}
\newcommand\bo{d{\overline \Omega}}
\newcommand\om{\omega}
\newcommand\SJ{{\mathbb H}_n\times {\mathbb C}^{(m,n)}}
\newcommand\DC{{\mathbb D}_n\times {\mathbb C}^{(m,n)}}
\newcommand\tr{\triangledown}
\newcommand\Hnm{{\mathbb H}_{n,m}}
\newcommand\Dnm{{\mathbb D}_{n,m}}
\newcommand\Hn{{\mathbb H}_n}
\newcommand\Cmn{{\mathbb C}^{(m,n)}}
\newcommand\ka{\kappa}

For brevity, we write $\Hnm:=\SJ.$ For a coordinate
$(\Om,Z)\in\Hnm$ with $\Om=(\om_{\mu\nu})\in {\mathbb H}_n$ and
$Z=(z_{kl})\in \Cmn,$ we put
\begin{eqnarray*}
\Om\,&=&\,X\,+\,iY,\quad\ \ X\,=\,(x_{\mu\nu}),\quad\ \
Y\,=\,(y_{\mu\nu})
\ \ \text{real},\\
Z\,&=&U\,+\,iV,\quad\ \ U\,=\,(u_{kl}),\quad\ \ V\,=\,(v_{kl})\ \
\text{real},\\
d\Om\,&=&\,(d\om_{\mu\nu}),\quad\ \ d{\overline
\Omega}=(d{\overline\omega}_{\mu\nu}), \\
dZ\,&=&\,(dz_{kl}),\quad\ \ d{\overline Z}=(d{\overline z}_{kl}),
\end{eqnarray*}

\begin{eqnarray*}
\PO\,=\,\left(\, { {1+\delta_{\mu\nu}} \over 2}\, {
{\partial}\over {\partial \om_{\mu\nu}} } \,\right),\quad
\POB\,=\,\left(\, { {1+\delta_{\mu\nu}}\over 2} \, {
{\partial}\over {\partial {\overline \om}_{\mu\nu} }  } \,\right),
\end{eqnarray*}

$$\PZ=\begin{pmatrix} {\partial}\over{\partial z_{11}} & \hdots &
 {\partial}\over{\partial z_{m1}} \\
\vdots&\ddots&\vdots\\
 {\partial}\over{\partial z_{1n}} &\hdots & {\partial}\over
{\partial z_{mn}} \end{pmatrix},\quad \PZB=\begin{pmatrix}
{\partial}\over{\partial {\overline z}_{11} }   &
\hdots&{ {\partial}\over{\partial {\overline z}_{m1} }  }\\
\vdots&\ddots&\vdots\\
{ {\partial}\over{\partial{\overline z}_{1n} }  }&\hdots &
 {\partial}\over{\partial{\overline z}_{mn} }  \end{pmatrix},$$

\noindent where $\delta_{ij}$ denotes the Kronecker delta symbol.

\newcommand\HC{{\mathbb H}\times{\mathbb C}}
\newcommand\ddx{{{\partial^2}\over{\partial x^2}}}
\newcommand\ddy{{{\partial^2}\over{\partial y^2}}}
\newcommand\ddu{{{\partial^2}\over{\partial u^2}}}
\newcommand\ddv{{{\partial^2}\over{\partial v^2}}}
\newcommand\px{{{\partial}\over{\partial x}}}
\newcommand\py{{{\partial}\over{\partial y}}}
\newcommand\pu{{{\partial}\over{\partial u}}}
\newcommand\pv{{{\partial}\over{\partial v}}}
\newcommand\pxu{{{\partial^2}\over{\partial x\partial u}}}
\newcommand\pyv{{{\partial^2}\over{\partial y\partial v}}}

\vskip 0.3cm C. L. Siegel [6] introduced the symplectic metric
$ds^2_n$ on ${\mathbb H}_n$ invariant under the action (1.1) of
$Sp(n,\BR)$ given by
\begin{eqnarray}
ds_n^2=\s \Big( Y^{-1}d\Om\,Y^{-1}d{\overline \Om} \Big)
\end{eqnarray}
and H. Maass [4] proved that the differential operator
\begin{eqnarray}
\Delta_n=\,4\,\s\left( Y \,\,^t\!\left( Y\POB\right) \PO\right)
%\s\left[Y\left\{\PU ^t\!\left(\PU\right)\right\}\right]
\end{eqnarray}
is the Laplacian of ${\mathbb H}_n$ for the symplectic metric
$ds_n^2.$ Here $\sigma(A)$ denotes the trace of a square matrix
$A$.

\vskip 0.3cm In [11], the author proved the following theorems.

\vskip 0.5cm \noindent {\bf Theorem\ A.} $ \textit{For any two
positive real numbers $A$ and $B$, the following metric}$
\begin{eqnarray}
ds_{n,m;A,B}^2&=&\,A\, \s\Big( Y^{-1}d\Om\,Y^{-1}d{\overline
\Om}\Big) \nonumber \\ && \ \ + \,B\,\bigg\{ \s\Big(
Y^{-1}\,^tV\,V\,Y^{-1}d\Om\,Y^{-1} \bo \Big)
 +\,\s\Big( Y^{-1}\,^t(dZ)\,\bz\Big)\\
&&\quad\quad -\s\Big( V\,Y^{-1}d\Om\,Y^{-1}\,^t(\bz)\Big)\,
-\,\s\Big( V\,Y^{-1}\bo\, Y^{-1}\,^t(dZ)\,\Big) \bigg\} \nonumber
\end{eqnarray}

\noindent $ \textit{is a Riemannian metric on $\Hnm$ which is
invariant under the action (1.2) of the Jacobi}$ \\
\noindent $ \textit{group $G^J$.}$

\vskip 0.3cm \noindent {\bf Theorem\ B.} $ \textit{For any two
positive real numbers $A$ and $B$, the Laplacian
$\Delta_{n,m;A,B}$ of}$
\\ \noindent $ \textit{$(\Hnm,ds^2_{n,m;A,B})$ is given by}$
\begin{eqnarray}
\Delta_{n,m;A,B}\,&=& \frac4A \,\bigg\{ \s\left(\,Y\,\,
^t\!\left(Y\POB\right)\PO\,\right)\, +\,\s\left(\,VY^{-1}\,^tV\,\,^t\!\left(Y\PZB\right)\,\PZ\,\right)\nonumber\\
& &\ \
+\,\s\left(V\,\,^t\!\left(Y\POB\right)\PZ\,\right)+\,\s\left(\,^tV\,\,^t\!\left(Y\PZB\right)\PO\,\right)\bigg\}\\
& & \ +\frac4B\,\,\s\left(\, Y\,\PZ\,\,{}^t\!\left(
\PZB\right)\,\right)\nonumber
\end{eqnarray}

\newcommand\w{\wedge}
\newcommand\OW{\overline{W}}
\newcommand\OP{\overline{P}}
\newcommand\OQ{\overline{Q}}
\newcommand\OZ{\overline{Z}}
\newcommand\Dn{{\mathbb D}_n}
\newcommand\la{\lambda}

\vskip 0.2cm Let \begin{equation*} G_*=SU(n,n)\cap
Sp(n,\BC)\end{equation*} be the symplectic group and let
\begin{equation*}
\BD_n=\left\{ W\in \BC^{(n,n)}\,|\ W=\,{}^tW,\ I_n-\OW W >
0\,\right\}
\end{equation*}
be the generalized unit disk. Then $G_*$ acts on $\Dn$
transitively by
$$\begin{pmatrix} P & Q\\
{\overline Q} & {\overline P}
\end{pmatrix}\cdot W=(PW+Q)(\OQ W+\OP)^{-1},$$
\noindent where $\begin{pmatrix} P & Q\\
{\overline Q} & {\overline P}
\end{pmatrix}\in G_*$ and $W\in\Dn.$
Using the Cayley transform of $\Dn$ onto $\Hn$, we can see
(cf.\,[6]) that
\begin{equation}
ds_*^2=4 \,\s \Big( (I_n-W{\overline W})^{-1}dW\,(I_n-\OW
W)^{-1}d\OW\,\Big)\end{equation} is a $G_*$-invariant K{\"a}hler
metric on $\BD_n$ and H. Maass [4] showed that its Laplacian is
given by
\begin{equation}
\Delta_*=\,\s \left( (I_n-W\OW)\,{ }^t\!\left( (I_n-W\OW)
\PWB\right)\PW\right).\end{equation} \vskip 0.2cm Let
\begin{equation*}
G^J_*=\left\{ \left( \begin{pmatrix} P & Q\\ \OQ & \OP
\end{pmatrix}, (\xi,{\bar \xi}\,;i\kappa)\right)\,\Big|
\ \begin{pmatrix} P & Q\\ \OQ & \OP
\end{pmatrix}
\in G_*,\ \xi\in \BC^{(m,n)},\ \kappa\in \BR^{(m,m)}\,\right\}
\end{equation*}
be the Jacobi group with the following multiplication
\begin{eqnarray*}
& &\,\ \   \left( \begin{pmatrix} P & Q\\ \OQ & \OP
\end{pmatrix}, (\xi,{\bar\xi}\,;i\kappa)\right) \cdot \left( \begin{pmatrix} P' & Q'\\
{\overline {Q'}} & {\overline {P'}}
\end{pmatrix}, (\xi',{\bar{\xi'}};i\kappa')\right)\\
&=& \left( \begin{pmatrix} P & Q\\ \OQ & \OP
\end{pmatrix}\,\begin{pmatrix} P' & Q'\\ {\overline{Q'}} &
{\overline{P'}}
\end{pmatrix},\,({\tilde \xi}+\xi',{\tilde
{\bar{\xi'}}}+{\overline{\xi'}};i\kappa+i\kappa'+{\tilde
\xi}\,{}^t{\overline{\xi'}}-{\tilde{\bar{\xi'}}} \,{}^t\xi')
\right),
\end{eqnarray*}
where ${\tilde\xi}=\xi P'+{\bar\xi}\,{\overline{Q'}} $ and
${\tilde {\bar \xi}} =\xi Q'+{\bar\xi}\,{\overline{P'}}.$ Then we
have the {\it canonical action} of $G_*^J$ on the $
\textit{Siegel}$-$ \textit{Jacobi disk}$ $\Dn\times \BC^{(m,n)}$
(see (2.6)) given by
\begin{equation}
\left(\begin{pmatrix} P & Q\\
{\overline Q} & {\overline P}
\end{pmatrix},\left( \xi, {\bar{\xi}};\,i\kappa\right)\right)\cdot
(W,\eta)=\Big((PW+Q)(\OQ W+\OP)^{-1},(\eta+\xi W+{\bar{\xi}})(\OQ
W+\OP)^{-1}\Big),\end{equation}

\noindent where $W\in\BD_n$ and $\eta\in\BC^{(m,n)}.$

\vskip 0.2cm For brevity, we write $\Dnm:=\DC.$ For a coordinate
$(W,\eta)\in\Dnm$ with $W=(w_{\mu\nu})\in {\mathbb D}_n$ and
$\eta=(\eta_{kl})\in \Cmn,$ we put
\begin{eqnarray*}
dW\,&=&\,(dw_{\mu\nu}),\quad\ \ d{\overline W}\,=\,(d{\overline w}_{\mu\nu}),\\
d\eta\,&=&\,(d\eta_{kl}),\quad\ \
d{\overline\eta}\,=\,(d{\overline\eta}_{kl})
\end{eqnarray*}
and
\begin{eqnarray*}
\PW\,=\,\left(\, { {1+\delta_{\mu\nu}} \over 2}\, {
{\partial}\over {\partial w_{\mu\nu}} } \,\right),\quad
\PWB\,=\,\left(\, { {1+\delta_{\mu\nu}}\over 2} \, {
{\partial}\over {\partial {\overline w}_{\mu\nu} }  } \,\right),
\end{eqnarray*}

$$\PE=\begin{pmatrix} {\partial}\over{\partial \eta_{11}} & \hdots &
 {\partial}\over{\partial \eta_{m1}} \\
\vdots&\ddots&\vdots\\
 {\partial}\over{\partial \eta_{1n}} &\hdots & {\partial}\over
{\partial \eta_{mn}} \end{pmatrix},\quad \PEB=\begin{pmatrix}
{\partial}\over{\partial {\overline \eta}_{11} }   &
\hdots&{ {\partial}\over{\partial {\overline \eta}_{m1} }  }\\
\vdots&\ddots&\vdots\\
{ {\partial}\over{\partial{\overline \eta}_{1n} }  }&\hdots &
 {\partial}\over{\partial{\overline \eta}_{mn} }  \end{pmatrix}.$$

\vskip 0.2cm In this paper, we find the $G_*^J$-invariant
Riemannian metrics on $\Dnm$ and their Laplacians. In fact, we
prove the following theorems.

\newcommand\ot{\overline\eta}
\vskip 0.5cm
\begin{theorem} For any two positive real numbers $A$ and $B$,
the following metric $d{\tilde s}^2_{n,m;A,B}$ defined by
\begin{eqnarray*}
d{\tilde s}^2_{n,m;A,B}&=&4\,A\,\s\Big( (I_n-W\OW)^{-1}dW(I_n-\OW W)^{-1}\bw\,\Big) \hskip 1cm\\
& &\,+\,4\,B\,\bigg\{ \s\Big(
(I_n-W\OW)^{-1}\,{}^t(d\eta)\,\be\,\Big)\\
& & \quad\quad\quad
\,+\,\s\Big(  (\eta\OW-{\overline\eta})(I_n-W\OW)^{-1}dW(I_n-\OW W)^{-1}\,{}^t(d\ot)\Big)\\
& & \quad\quad\quad  +\,\s\Big( (\ot W-\eta)(I_n-\OW
W)^{-1}d\OW(I_n-W\OW)^{-1}\,{}^t(d\eta)\,\Big)    \\
& &\quad\quad\quad -\,\s\Big( (I_n-W\OW)^{-1}\,{}^t\eta\,\eta\,
(I_n-\OW W)^{-1}\OW
dW (I_n-\OW W)^{-1}d\OW \, \Big)\\
& &\quad\quad\quad -\,\s\Big( W(I_n-\OW W)^{-1}\,{}^t\ot\,\ot\,
(I_n-W\OW )^{-1}
dW (I_n-\OW W)^{-1}d\OW \,\Big)\\
& &\quad\quad\quad +\,\s\Big( (I_n-W\OW)^{-1}{}^t\eta\,\ot \,(I_n-W\OW)^{-1} dW (I_n-\OW W)^{-1} d\OW\,\Big)\\
& &\quad\quad\quad +\,\s\Big( (I_n-\OW)^{-1}\,{}^t\ot\,\eta\,\OW\,(I_n-W\OW)^{-1} dW (I_n-\OW W)^{-1} d\OW\,\Big)\\
& &\quad\quad\quad +\,\s \Big( (I_n-\OW)^{-1}(I_n-W)(I_n-\OW
W)^{-1}\,{}^t\ot\,\eta\,(I_n-\OW W)^{-1}\\
& &\qquad\qquad\quad\quad \times\, (I_n-\OW)(I_n-W)^{-1}dW
(I_n-\OW W)^{-1}d\OW\,\Big)\\
& &\quad\quad\quad -\,\s
\Big( (I_n-W\OW)^{-1}(I_n-W)(I_n-\OW)^{-1}\,{}^t\ot\,\eta\,(I_n-W)^{-1}\\
& & \qquad\qquad\quad\quad \times\,dW (I_n-\OW
W)^{-1}d\OW\,\Big)\bigg\}
\end{eqnarray*}
\noindent is a Riemannian metric on $\Dnm$ which is invariant
under the action (1.9) of the Jacobi group $G^J_*$. \vskip 0.2cm
We note that if $n=m=1$ and $A=B=1,$ we get
\begin{eqnarray*}
{\frac 14}\,d{\tilde s}_{1,1;1,1}^2&=& { {dW\,d\OW}\over
{(1-|W|^2)^2}}\,+\,{ 1 \over {(1-|W|^2)} }\,d\eta\,d\ot\\
& &+{ {(1+|W|^2)|\eta|^2-\OW \eta^2-W\ot^2}\over {(1-|W|^2)^3} }\,dW\,d\OW\\
& & + { {\eta\OW -\ot}\over {(1-|W|^2)^2} }\,dWd\ot\,+\,{ {\ot W
-\eta}\over {(1-|W|^2)^2} }\,d\OW d\eta.
\end{eqnarray*}
\end{theorem}

%\newpage
\vskip 0.3cm
\begin{theorem} For any two positive real numbers $A$ and $B$,
the Laplacian ${\tilde\Delta}_{n,m;A,B}$ of $(\Dnm,d{\tilde
s}_{n,m;A,B}^2)$ is given by
\begin{eqnarray*}
{\tilde\Delta}_{n,m;A,B}\,&=&\,\frac1A\,\bigg\{ \s\left(
(I_n-W\OW)\,{}^t\!\left(
(I_n-W\OW)\PWB\right)\PW\right)\,\\
& & \quad \quad +\,\s \left(\,{}^t(\eta-\ot\,W)\,{}^t\!\left(
\PEB\right)
(I_n-\OW W)\PW  \right)\,\\
& & \quad \quad +\,\s\left( (\ot-\eta\,\OW)\,{}^t\!\left(
(I_n-W\OW)\PWB\right)\PE\right)\\
& &\quad \quad -\,\s\left( \eta \OW
(I_n-W\OW)^{-1}\,{}^t\eta\,\,{}^t\!\left(\PEB\right)(I_n-\OW
W)\PE\right)\\
& &\quad \quad -\,\s\left( \ot W (I_n-\OW W)^{-1}
\,{}^t\ot\,\,{}^t\!\left(\PEB\right)(I_n-\OW
W)\PE\right)\\
& &\quad \quad +\,\s\left( \ot
(I_n-W\OW)^{-1}{}^t\eta\,{}^t\!\left( \PEB\right)
(I_n-\OW W)\PE \right)\\
& & \quad \quad +\,\s\left( \eta\,\OW W (I_n-\OW
W)^{-1}\,{}^t\ot\,\,{}^t\!\left( \PEB\right) (I_n-\OW W)\PE
\right)\bigg\}\\
& & \  +\,\frac1B \,\s\left( (I_n-\OW W) \PE\,{}^t\!\left(
\PEB\right)\right).
\end{eqnarray*}
\newcommand\ddww{{{\partial^2}\over{\partial W\partial \OW}}}
\newcommand\ddtt{{{\partial^2}\over{\partial \eta\partial \ot}}}
\newcommand\ddwe{{{\partial^2}\over{\partial W\partial \ot}}}
\newcommand\ddew{{{\partial^2}\over{\partial \OW\partial \eta}}}

\vskip 0.2cm We note that if $n=m=1$ and $A=B=1,$ we get
\begin{eqnarray*}
{\tilde\Delta}_{1,1;1,1}&=&\ \ \Big(1-|W|^2\Big)^2 \ddww\,+\,\Big(1-|W|^2\Big)\,\ddtt\\
& &+\,\Big(1-|W|^2)(\eta-\ot\,W\Big)\,\ddwe\,+\,\Big(1-|W|^2)(\ot-\eta\,\OW\Big)\,\ddew\\
&
&-\Big(\OW\,\eta^2+W\ot^2\Big)\,\ddtt\,+\,\Big(1+|W|^2\Big)|\eta|^2\,\ddtt.
\end{eqnarray*}
\end{theorem}

\vskip 0.2cm The main ingredients for the proof of Theorem 1.1 and
Theorem 1.2 are the partial Cayley transform, Theorem A and
Theorem B. The paper is organized as follows. In Section 2, we
review the partial Cayley transform that was dealt with in [12].
In Section 3, we prove Theorem 1.1. In Section 4, we prove Theorem
1.2. In the final section we briefly remark the theory of harmonic
analysis
%spectral theory of the Laplacian ${\tilde \Delta}_{n,m;A,B}$
on the Siegel-Jacobi disk.

\vskip 0.2cm  \noindent {\textsc{Notations}\,:} \ We denote by
$\BR$ and $\BC$ the field of real numbers, and the field of
complex numbers respectively. The symbol ``:='' means that the
expression on the right is the definition of that on the left. For
two positive integers $k$ and $l$, $F^{(k,l)}$ denotes the set of
all $k\times l$ matrices with entries in a commutative ring $F$.
For a square matrix $A\in F^{(k,k)}$ of degree $k$, $\sigma(A)$
denotes the trace of $A$. For $\Omega\in {\mathbb H}_g,\
\text{Re}\,\Omega$ ({\it resp.}\ $\textrm{Im}\,\Omega)$ denotes
the real ({\it resp.}\ imaginary) part of $\Omega.$ For any $M\in
F^{(k,l)},\ ^t\!M$ denotes the transpose matrix of $M$. For a
matrix $A\in F^{(k,k)}$ and $B\in F^{(k,l)},$ we write
$A[B]=\,^tBAB$. $I_n$ denotes the identity matrix of degree $n$.

\vskip 0.5cm
\begin{section}{{\bf \ The partial Cayley transform}}
\setcounter{equation}{0}
\end{section}
\vskip 0.2cm In this section, we review the partial Cayley
transform\,[12] of $\Dnm$ onto $\Hnm$ needed for the proof of
Theorem 1.1 and Theorem 1.2.

\vskip 0.2cm  We can identify an element $g=(M,(\la,\mu;\kappa))$
of $G^J,\ M=\begin{pmatrix} A&B\\
C&D\end{pmatrix}\in Sp(n,\BR)$ with the element
\begin{equation*}
\begin{pmatrix} A & 0 & B & A\,^t\mu-B\,^t\la  \\ \la & I_m & \mu
& \kappa \\ C & 0 & D & C\,^t\mu-D\,^t\la \\ 0 & 0 & 0 & I_m
\end{pmatrix}
\end{equation*}
of $Sp(m+n,\BR).$ \vskip 0.3cm We set
\begin{equation*}
T_*={1\over {\sqrt 2}}\,
\begin{pmatrix} I_{m+n} & I_{m+n}\\ iI_{m+n} & -iI_{m+n}
\end{pmatrix}.
\end{equation*}
We now consider the group $G_*^J$ defined by
\begin{equation*}
G_*^J:=T_*^{-1}G^JT_*.
\end{equation*}
If $g=(M,(\la,\mu;\kappa))\in G^J$ with $M=\begin{pmatrix} A&B\\
C&D\end{pmatrix}\in Sp(n,\BR)$, then $T_*^{-1}gT_*$ is given by
\begin{equation}
T_*^{-1}gT_*=
\begin{pmatrix} P_* & Q_*\\ {\overline Q}_* & {\overline P}_*
\end{pmatrix},
\end{equation}
where
\begin{equation*}
P_*=
\begin{pmatrix} P & {\frac 12} \left\{ Q\,\,{}^t(\la+i\mu)-P\,\,{}^t(\la-i\mu)\right\}\\
{\frac 12} (\la+i\mu) & I_h+i{\frac \kappa 2}
\end{pmatrix},
\end{equation*}

\begin{equation*}
Q_*=
\begin{pmatrix} Q & {\frac 12} \left\{ P\,\,{}^t(\la-i\mu)-Q\,\,{}^t(\la+i\mu)\right\}\\
{\frac 12} (\la-i\mu) & -i{\frac \kappa 2}
\end{pmatrix},
\end{equation*}
and $P,\,Q$ are given by the formulas
\begin{equation}
P= {\frac 12}\,\left\{ (A+D)+\,i\,(B-C)\right\}
\end{equation}
and
\begin{equation}
 Q={\frac
12}\,\left\{ (A-D)-\,i\,(B+C)\right\}.
\end{equation}
From now on, we write
\begin{equation*}
\left(\begin{pmatrix} P & Q\\ {\overline Q} & {\overline P}
\end{pmatrix},\left( {\frac 12}(\la+i\mu),\,{\frac 12}(\la-i\mu);\,-i{\kappa\over 2}\right)\right):=
\begin{pmatrix} P_* & Q_*\\ {\overline Q}_* & {\overline P}_*
\end{pmatrix}.
\end{equation*}
In other words, we have the relation
\begin{equation*}
T_*^{-1}\left( \begin{pmatrix} A&B\\
C&D\end{pmatrix},(\la,\mu;\kappa)
\right)T_*=  \left(\begin{pmatrix} P & Q\\
{\overline Q} & {\overline P}
\end{pmatrix},\left(
{\frac 12}(\la+i\mu),\,{\frac 12}(\la-i\mu);\,-i{\kappa\over 2}
\right)\right).
\end{equation*}
Let
\begin{equation*}
H_{\BC}^{(n,m)}:=\left\{ (\xi,\eta\,;\zeta)\,|\
\xi,\eta\in\BC^{(m,n)},\ \zeta\in\BC^{(m,m)},\
\zeta+\eta\,{}^t\xi\ \textrm{symmetric}\,\right\}
\end{equation*}
be the complex Heisenberg group endowed with the following
multiplication
\begin{equation*}
(\xi,\eta\,;\zeta)\circ
(\xi',\eta';\zeta'):=(\xi+\xi',\eta+\eta'\,;\zeta+\zeta'+
\xi\,{}^t\eta'-\eta\,{}^t\xi')).
\end{equation*}
We define the semidirect product
\begin{equation*}
%SLH(2g,\BC):=
SL(2n,\BC)\ltimes H_{\BC}^{(n,m)}
\end{equation*}
endowed with the following multiplication
\begin{eqnarray*}
& & \left( \begin{pmatrix} P & Q\\ R & S
\end{pmatrix}, (\xi,\eta\,;\zeta)\right)\cdot \left( \begin{pmatrix} P' & Q'\\
R' & S'
\end{pmatrix}, (\xi',\eta';\zeta')\right)\\
&=& \left( \begin{pmatrix} P & Q\\ R & S
\end{pmatrix}\,\begin{pmatrix} P' & Q'\\ R' & S'
\end{pmatrix},\,({\tilde \xi}+\xi',{\tilde
\eta}+\eta';\zeta+\zeta'+{\tilde \xi}\,{}^t\eta'-{\tilde
\eta}\,{}^t\xi')  \right),
\end{eqnarray*}
where ${\tilde\xi}=\xi P'+\eta R'$ and ${\tilde \eta}=\xi Q'+\eta
S'.$

\vskip 0.2cm If we identify $H_{\BR}^{(n,m)}$ with the subgroup
$$\left\{ (\xi,{\overline \xi};i\kappa)\,|\ \xi\in\BC^{(m,n)},\
\ka\in\BR^{(m,m)}\,\right\}$$ of $H_{\BC}^{(n,m)},$ we have the
following inclusion
$$G_*^J\subset SU(n,n)\ltimes H_{\BR}^{(n,m)}\subset SL(2n,\BC)\ltimes
H_{\BC}^{(n,m)}.$$ We define the mapping $\Theta:G^J\lrt G_*^J$ by
\begin{equation}\Theta
\left( \begin{pmatrix} A&B\\
C&D\end{pmatrix},(\la,\mu;\kappa) \right):=\left(\begin{pmatrix} P
& Q\\ {\overline Q} & {\overline P}
\end{pmatrix},\left(
{\frac 12}(\la+i\mu),\,{\frac 12}(\la-i\mu);\,-i{\kappa\over 2}
\right)\right),\end{equation} where $P$ and $Q$ are given by (2.2)
and (2.3). We can see that if $g_1,g_2\in G^J$, then
$\Theta(g_1g_2)=\Theta(g_1)\Theta(g_2).$

\vskip 0.2cm According to [10], p.\,250, $G_*^J$ is of the
Harish-Chandra type\,(cf.\,[5],\,p.\,118). Let
$$g_*=\left(\begin{pmatrix} P & Q\\
{\overline Q} & {\overline P}
\end{pmatrix},\left( \la, \mu;\,\kappa\right)\right)$$
be an element of $G_*^J.$ Since the Harish-Chandra decomposition
of an element $\begin{pmatrix} P & Q\\ R & S
\end{pmatrix}$ in $SU(n,n)$ is given by
\begin{equation*}
\begin{pmatrix} P & Q\\ R & S
\end{pmatrix}=\begin{pmatrix} I_n & QS^{-1}\\ 0 & I_n
\end{pmatrix} \begin{pmatrix} P-QS^{-1}R & 0\\ 0 & S
\end{pmatrix} \begin{pmatrix} I_n & 0\\ S^{-1}R & I_n
\end{pmatrix},
\end{equation*}
the $P_*^+$-component of the following element
$$g_*\cdot\left( \begin{pmatrix} I_n & W\\ 0 & I_n
\end{pmatrix}, (0,\eta;0)\right),\quad W\in \BD_n$$
of $SL(2n,\BC)\ltimes H_{\BC}^{(n,m)}$ is given by
\begin{equation}
\left( \begin{pmatrix} I_n & (PW+Q)(\OQ W+\OP)^{-1}
\\ 0 & I_n
\end{pmatrix},\,\left(0,\,(\eta+\la W+\mu)(\OQ W+\OP)^{-1}\,;0\right)\right).
\end{equation}

\vskip 0.2cm We can identify $\Dnm$ with the subset
\begin{equation*}
\left\{ \left( \begin{pmatrix} I_n & W\\ 0 & I_n
\end{pmatrix}, (0,\eta;0)\right)\,\Big|\ W\in\BD_n,\
\eta\in\BC^{(m,n)}\,\right\}\end{equation*} of the
complexification of $G_*^J.$ Indeed, $\Dnm$ is embedded into
$P_*^+$ given by
\begin{equation*}
P_*^+=\left\{\,\left( \begin{pmatrix} I_n & W\\ 0 & I_n
\end{pmatrix}, (0,\eta;0)\right)\,\Big|\ W=\,{}^tW\in \BC^{(n,n)},\
\eta\in\BC^{(m,n)}\ \right\}.
\end{equation*}
This is a generalization of the Harish-Chandra
embedding\,(cf.\,[5],\,p.\,119). Then we get the {\it canonical
transitive action} of $G_*^J$ on $\Dnm$ defined by
\begin{equation}
\left(\begin{pmatrix} P & Q\\
{\overline Q} & {\overline P}
\end{pmatrix},\left( \xi, {\overline\xi};\,i\kappa\right)\right)\cdot
(W,\eta)=\Big((PW+Q)(\OQ W+\OP)^{-1},(\eta+\xi
W+{\overline\xi})(\OQ W+\OP)^{-1}\Big),\end{equation}

\noindent where $\begin{pmatrix} P & Q\\
{\overline Q} & {\overline P}
\end{pmatrix}\in G_*,\ \xi\in \BC^{(m,n)},\ \k\in\BR^{(m,m)}$ and
$(W,\eta)\in\Dnm.$

\vskip 0.2cm The author [12] proved that the action (1.2) of $G^J$
on $\Hnm$ is compatible with the action (2.6) of $G_*^J$ on $\Dnm$
through the {\it partial Cayley transform} $\Phi:\BD_{n,m}\lrt
\BH_{n,m}$ defined by
\begin{equation}
\Phi(W,\eta):=\Big(
i(I_n+W)(I_n-W)^{-1},\,2\,i\,\eta\,(I_n-W)^{-1}\Big).
\end{equation}
In other words, if $g_0\in G^J$ and $(W,\eta)\in\BD_{n,m}$,
\begin{equation}
g_0\cdot\Phi(W,\eta)=\Phi(g_*\cdot (W,\eta)),
\end{equation}
where $g_*=T_*^{-1}g_0 T_*$. $\Phi$ is a biholomorphic mapping of
$\Dnm$ onto $\Hnm$ which gives the partially bounded realization
of $\Hnm$ by $\Dnm$. The inverse of $\Phi$ is
\begin{equation*}
\Phi^{-1}(\Omega,Z)=\Big(
(\Omega-iI_n)(\Omega+iI_n)^{-1},\,Z(\Omega+iI_n)^{-1}\Big).
\end{equation*}

\vskip 0.5cm
\begin{section}{{\bf \ Proof\ of\ Theorem\ 1.1}}
\setcounter{equation}{0}
\end{section}
\vskip 0.1cm For $(W,\eta)\in \Dnm,$ we write
\begin{equation*}
(\Omega,Z):=\Phi(W,\eta).
\end{equation*}
Thus
\begin{equation}
\Om=i(I_n+W)(I_n-W)^{-1},\qquad Z=2\,i\,\eta\,(I_n-W)^{-1}.
\end{equation}

Since
\begin{equation*}
d(I_n-W)^{-1}=(I_n-W)^{-1}dW\,(I_n-W)^{-1}
\end{equation*}
and
\begin{equation*}
I_n+(I_n+W)(I_n-W)^{-1}=2\,(I_n-W)^{-1},
\end{equation*}
we get the following formulas from (3.1)

\begin{eqnarray}
Y&=&{1 \over {2\,i}}\,(\Om-{\overline
\Om}\,)=(I_n-W)^{-1}(I_n-W\OW\,)(I_n-\OW\,)^{-1},\\
V&=&{1 \over {2\,i}}\,(Z-\OZ\,)=\eta\, (I_n-W)^{-1}+\ot\,
(I_n-\OW\,)^{-1},\\
d\Om &=& 2\,i\,(I_n-W)^{-1}dW\,(I_n-W)^{-1},\\
dZ&=&2\,i\,\Big\{ d\eta+\eta\,(I_n-W)^{-1}dW\,\Big\}(I_n-W)^{-1}.
\end{eqnarray}

\vskip 0.2cm According to the formulas (3.2) and (3.4), we obtain
\begin{equation}
\s\Big( Y^{-1}d\Om \,Y^{-1}d{\overline\Om}\,\Big)=4\,\s\Big(
(I_n-W\OW\,)^{-1}dW (I_n-\OW W)^{-1}d\OW\,\Big).
\end{equation}

\vskip 0.2cm From the formulas (3.2)-(3.4), we get
\begin{equation}
\s\Big( Y^{-1}\,{}^tV VY^{-1}d\Om
\,Y^{-1}d{\overline\Om}\,\Big)=(a)+(b)+(c)+(d),
\end{equation}
where
\begin{eqnarray*}
(a):&=& \,4\,\s \Big( (I_n-W\OW\,)^{-1}\,{}^t\eta\,\eta\,(I_n-\OW
W)^{-1}(I_n-\OW\,)(I_n-W)^{-1}\\
& & \qquad\quad  \times\, dW\,(I_n-\OW W)^{-1}d\OW\,\Big),\\
(b):&=&\,4\,\s\Big(
(I_n-W\OW\,)^{-1}\,{}^t\eta\,\ot\,(I_n-W\OW\,)^{-1}dW\,(I_n-\OW
W)^{-1}d\OW\,\Big),\\
 (c):&=&\,4\,\s \Big(
(I_n-\OW\,)^{-1}(I_n-W)(I_n-\OW
W)^{-1}\,{}^t\ot \,\eta\,(I_n-\OW W)^{-1}\\
& &\qquad\quad\times (I_n-\OW\,)(I_n-W)^{-1}dW\,(I_n-\OW
W)^{-1}d\OW\,\Big),\\
(d):&=&\,4\,\s \Big( (I_n-\OW\,)^{-1}(I_n-W)(I_n-\OW
W)^{-1}\,{}^t\ot \,\ot \,(I_n-W\OW\,)^{-1}\\
& &\qquad\quad\times\, dW\,(I_n-\OW W)^{-1}d\OW\,\Big).
\end{eqnarray*}

\vskip 0.2cm
According to the formulas (3.2) and (3.5), we get
\begin{equation}
\s\Big( Y^{-1}\,{}^t(dZ)\,d\OZ\,\Big)=(e)+(f)+(g)=(h),
\end{equation}
where
\begin{eqnarray*}
(e):&=& \,4\,\s \Big(
(I_n-W\OW\,)^{-1}\,{}^t(d\eta)\,d\ot\,\Big),\\
(f):&=& \,4\,\s \Big(
(I_n-W\OW\,)^{-1}\,dW\,(I_n-W)^{-1}\,{}^t\eta\,d\ot\,\Big),\\
(g):&=&\,4\,\s \Big(
(I_n-W\OW\,)^{-1}\,{}^t(d\eta)\,\ot\,(I_n-\OW\,)^{-1}d\OW\,\Big),\\
(h):&=&\,4\,\s \Big(
(I_n-W\OW\,)^{-1}\,dW\,(I_n-W)^{-1}\,{}^t\eta\,\ot\,(I_n-\OW\,)^{-1}d\OW\,\Big).
\end{eqnarray*}

\vskip 0.2cm
From the formulas (3.2)-(3.5), we get
\begin{equation}
-\s\Big( V Y^{-1} d\Om \,Y^{-1}
\,{}^t(d\OZ\,)\,\Big)=(i)+(j)+(k)+(l),
\end{equation}
where
\begin{eqnarray*}
(i):&=& \,-4\,\s \Big(
\eta\,(I_n-W)^{-1}(I_n-\OW\,)(I_n-W\OW\,)^{-1}dW\,(I_n-\OW
W)^{-1}\,{}^t(d\ot\,)\Big),\\
(j):&=&\,-4\,\s \Big( (I_n-\OW\,)^{-1}\,{}^t\ot\,\eta\,(I_n-W)^{-1}(I_n-\OW\,)(I_n-W\OW\,)^{-1}\\
& &\qquad\quad\times\, dW\,(I_n-\OW W)^{-1}d\OW\,\Big),\\
 (k):&=&\,-4\,\s \Big(\,\ot\,(I_n-W\OW\,)^{-1}dW\,(I_n-\OW
W)^{-1}\,{}^t(d\ot\,)\Big),
\\
(l):&=&\,-4\,\s \Big(
(I_n-\OW\,)^{-1}\,{}^t\ot\,\ot\,(I_n-W\OW\,)^{-1}  dW\,(I_n-\OW
W)^{-1}d\OW\,\Big).
\end{eqnarray*}

\vskip 0.2cm Conjugating the formula (3.9), we get
\begin{equation}
-\s\Big( V Y^{-1} d{\overline\Om} \,Y^{-1}
\,{}^t(dZ)\,\Big)=(m)+(n)+(o)+(p),
\end{equation}
where
\begin{eqnarray*}
(m):&=& \,-4\,\s \Big( \ot\,(I_n-\OW\,)^{-1}(I_n-W)(I_n-\OW
W)^{-1}d\OW\,(I_n-W\OW\,)^{-1}\,{}^t(d\eta\,)\Big),\\
(n):&=&\,-4\,\s \Big( (I_n-W)^{-1}\,{}^t\eta\,\ot\,(I_n-\OW\,)^{-1}(I_n-W)(I_n-\OW W)^{-1}\\
& &\qquad\quad\times\, d\OW\,(I_n-W\OW\,)^{-1}dW\,\Big),\\
 (o):&=&\,-4\,\s \Big(\,\eta\,(I_n-\OW W)^{-1}d\OW\,(I_n-W\OW\,)^{-1}\,{}^t(d\eta\,)\Big),
\\
(p):&=&\,-4\,\s \Big( (I_n-W)^{-1}\,{}^t\eta\,\eta\,(I_n-\OW
W)^{-1} d\OW\,(I_n-W\OW\,)^{-1}dW\,\Big).
\end{eqnarray*}

\vskip 0.2cm If we add the formulas $(f),\,(i)$ and $(k)$, we get

\begin{equation}
(f)+(i)+(k)=\,4\,\s\Big( (\eta
\OW-\ot\,)(I_n-W\OW\,)^{-1}dW\,(I_n-\OW
W)^{-1}\,{}^t(d\ot\,)\Big).
\end{equation}

\vskip 0.2cm Indeed, transposing the matrix inside the formula
$(f)$, we can express the formula $(f)$ as

\begin{equation*}
(f)=\,4\,\s\Big( \eta\,(I_n-W)^{-1}\,dW (I_n-\OW
W)^{-1}\,{}^t(\be)\Big)
\end{equation*}

\vskip 0.2cm \noindent
and adding the formulas $(f)$ and $(i)$
together with the formula $(k)$, we get the formula (3.11) because
\begin{eqnarray*}
& & (I_n-W)^{-1}-(I_n-W)^{-1}(I_n-\OW\,)(I_n-W\OW\,)^{-1}\\
&=&(I_n-W)^{-1}\Big\{ (I_n-W\OW\,)-(I_n-\OW\,)\Big\}(I_n-W\OW\,)^{-1}\\
&=&(I_n-W)^{-1}(I_n-W)\OW\,(I_n-W\OW\,)^{-1}\\
&=& \OW\,(I_n-W\OW\,)^{-1}.
\end{eqnarray*}

\vskip 0.3cm If we add the formulas $(g),\,(m)$ and $(o)$, we get

\begin{equation}
(g)+(m)+(o)=\,4\,\s\Big( (\ot \,W-\eta\,)(I_n-\OW
W)^{-1}d\OW\,(I_n-W\OW\,)^{-1}\,{}^t(d\eta\,)\Big).
\end{equation}

\vskip 0.2cm Indeed, we can express the formula $(g)$ as

\begin{equation*}
(g)=\,4\,\s\Big(
\ot\,(I_n-\OW\,)^{-1}d\OW\,(I_n-W\OW\,)^{-1}\,{}^t(d\eta)\Big)
\end{equation*}

\vskip 0.2cm\noindent and adding the formulas $(g)$ and $(m)$
together with the formula $(o)$, we get the formula (3.12) because
\begin{eqnarray*}
& & (I_n-\OW\,)^{-1}-(I_n-\OW\,)^{-1}(I_n-W)(I_n-\OW\,W)^{-1}\\
&=&(I_n-\OW\,)^{-1}\Big\{ (I_n-\OW W)-(I_n-W)\Big\}(I_n-\OW W)^{-1}\\
&=&(I_n-\OW\,)^{-1}(I_n-\OW\,)W(I_n-\OW W)^{-1}\\
&=& W(I_n-\OW W)^{-1}.
\end{eqnarray*}

\vskip 0.3cm If we add the formulas $(a)$ and $(p)$, we get

\begin{eqnarray}
(a)+(p)&=&\,-4\,\s \Big( (I_n-W\OW\,)^{-1}\,{}^t\eta\,\eta
\,(I_n-\OW\,W)^{-1}\OW\\
& & \quad\qquad \times\, \,
dW\,(I_n-\OW\,W)^{-1}d\OW\,\Big).\nonumber
\end{eqnarray}

\vskip 0.2cm Indeed, transposing the matrix inside the formula
$(p)$, we can express the formula $(p)$ as

\begin{equation*}
(p)=\,-4\,\s\Big( (I_n-W\OW\,)^{-1}\,{}^t\eta\,
\eta\,(I_n-W)^{-1}dW \,(I_n-\OW W)^{-1}\,\bw\Big)
\end{equation*}

\vskip 0.2cm\noindent and adding the formulas $(a)$ and $(p)$, we
get the formula (3.13) because

\begin{eqnarray*}
& & (I_n-\OW\,W)^{-1}(I_n-\OW\,)(I_n-W)^{-1}-(I_n-W)^{-1}\\
&=&(I_n-\OW\,W)^{-1}\Big\{ (I_n-\OW )-(I_n-\OW\,W)\Big\}(I_n- W)^{-1}\\
&=&(I_n-\OW\,W)^{-1}(-\OW)(I_n-W)(I_n- W)^{-1}\\
&=& -(I_n-\OW W)^{-1}\OW.
\end{eqnarray*}

\vskip 0.3cm Adding the formulas $(d)$ and $(l)$, we get

\begin{eqnarray}
(d)+(l)&=&\,-4\,\s \Big( W(I_n-\OW W)^{-1}\,{}^t\ot\,\ot
\,(I_n-W\OW\,)^{-1}\\
& & \quad\qquad \times\, dW\,(I_n-\OW\,W)^{-1}d\OW\,\Big)\nonumber
\end{eqnarray}
\vskip 0.2cm \noindent because
\begin{eqnarray*}
& & (I_n-\OW\,)^{-1}(I_n-W)(I_n-\OW\,W)^{-1}-(I_n-\OW\,)^{-1}\\
&=&(I_n-\OW\,)^{-1}\Big\{ (I_n-W)-(I_n-\OW W)\Big\}(I_n-\OW W)^{-1}\\
&=&(I_n-\OW\,)^{-1}(I_n-\OW\,)(-W)(I_n-\OW W)^{-1}\\
&=& -W(I_n-\OW W)^{-1}.
\end{eqnarray*}

\vskip 0.3cm Adding the formulas $(h)$ and $(j)$, we get

\begin{eqnarray}
(h)+(j)&=&\,4\,\s \Big( (I_n-\OW \,)^{-1}\,{}^t\ot\,\eta
\,\OW\,(I_n-W\OW\,)^{-1}\\
& & \quad\qquad \times\,
dW\,(I_n-\OW\,W)^{-1}d\OW\,\Big).\nonumber
\end{eqnarray}

\vskip 0.2cm Indeed, transposing the matrix inside the formula
$(h)$, we can express the formula $(h)$ as

\begin{equation*}
(h)=\,4\,\s\Big( (I_n-\OW\,)^{-1}\,{}^t\ot\,\eta\,(I_n-W)^{-1}\,dW
(I_n-\OW W)^{-1}\,\bw\Big)
\end{equation*}

\vskip 0.2cm \noindent and adding the formulas $(h)$ and $(j)$, we
get the formula (3.15) because

\begin{eqnarray*}
& & (I_n-W)^{-1}-(I_n-W)^{-1}(I_n-\OW\,)(I_n-W\OW\,)^{-1}\\
&=&(I_n-W)^{-1}\Big\{ (I_n-W\OW\,)-(I_n-\OW\,)\Big\}(I_n-W\OW\,)^{-1}\\
&=&(I_n-W)^{-1}(I_n-W)\OW\,(I_n-W\OW\,)^{-1}\\
&=& \OW\,(I_n-W\OW\,)^{-1}.
\end{eqnarray*}

\vskip 0.3cm Transposing the matrix inside the formula $(n)$, we
get

\begin{eqnarray}
(n)&=&\,-4\,\s \Big( (I_n-W\OW
\,)^{-1}(I_n-W)(I_n-\OW\,)^{-1}\,{}^t\ot\,\eta
\,(I_n-W)^{-1}\\
& & \quad\qquad \times\,
dW\,(I_n-\OW\,W)^{-1}d\OW\,\Big).\nonumber
\end{eqnarray}

From the formulas (3.7)-(3.16), we obtain
\begin{eqnarray*}
&& \quad\s\Big( Y^{-1}\,^tV\,V\,Y^{-1}d\Om\,Y^{-1} \bo \Big)
 +\,\s\Big( Y^{-1}\,^t(dZ)\,\bz\Big)\\
&&\quad\ \ -\s\Big( V\,Y^{-1}d\Om\,Y^{-1}\,^t(\bz)\Big)\,
-\,\s\Big(
V\,Y^{-1}\bo\, Y^{-1}\,^t(dZ)\,\Big)\\
&&=(a)+(b)+(c)+(d)+\cdots+(m)+(n)+(o)+(p)
\end{eqnarray*}

\begin{eqnarray*}
&&=4\,\s\Big(
(I_n-W\OW)^{-1}\,{}^t(d\eta)\,\be\,\Big)\\
& & \quad
\,+\,4\,\s\Big(  (\eta\OW-{\overline\eta})(I_n-W\OW)^{-1}dW(I_n-\OW W)^{-1}\,{}^t(d\ot)\Big)\\
& & \quad  +\,4\,\s\Big( (\ot W-\eta)(I_n-\OW
W)^{-1}d\OW(I_n-W\OW)^{-1}\,{}^t(d\eta)\,\Big)    \\
& &\quad -\,4\,\s\Big( (I_n-W\OW)^{-1}\,{}^t\eta\,\eta\, (I_n-\OW
W)^{-1}\OW
dW (I_n-\OW W)^{-1}d\OW \, \Big)\\
& &\quad -\,4\,\s\Big( W(I_n-\OW W)^{-1}\,{}^t\ot\,\ot\, (I_n-W\OW
)^{-1}
dW (I_n-\OW W)^{-1}d\OW \,\Big)\\
& &\quad +\,4\,\s\Big( (I_n-W\OW)^{-1}{}^t\eta\,\ot \,(I_n-W\OW)^{-1} dW (I_n-\OW W)^{-1} d\OW\,\Big)\\
& &\quad +\,4\,\s\Big( (I_n-\OW)^{-1}\,{}^t\ot\,\eta\,\OW\,(I_n-W\OW)^{-1} dW (I_n-\OW W)^{-1} d\OW\,\Big)\\
& &\quad +\,4\,\s \Big( (I_n-\OW)^{-1}(I_n-W)(I_n-\OW
W)^{-1}\,{}^t\ot\,\eta\,(I_n-\OW W)^{-1}\\
& &\qquad\quad\quad \times\, (I_n-\OW)(I_n-W)^{-1}dW
(I_n-\OW W)^{-1}d\OW\,\Big)\\
& &\quad -\,4\,\s
\Big( (I_n-W\OW)^{-1}(I_n-W)(I_n-\OW)^{-1}\,{}^t\ot\,\eta\,(I_n-W)^{-1}\\
& & \qquad\quad\quad \times\,dW (I_n-\OW W)^{-1}d\OW\,\Big).
\end{eqnarray*}
\vskip 0.3cm Consequently the complete proof follows from the
above formula, the formula (3.6), Theorem A and the fact that the
action (1.2) of $G^J$ on $\Hnm$ is compatible with the action
(2.6) of $G_*^J$ on $\Dnm$ through the partial Cayley
transform.\hfill $\square$

\vskip 1.05cm
\begin{section}{{\bf \ Proof\ of\ Theorem\ 1.2}}
\setcounter{equation}{0}
\end{section}
\vskip 0.1cm From the formulas (3.1),\,(3.4) and (3.5), we get

\begin{equation}
\PO=\,{1\over {2\,i}}\,(I_n-W)\left[ \,{}^t\left\{
(I_n-W)\,\PW\right\} -\,{}^t\left\{ \,{}^t\eta\,\,{}^t\left(
\PE\right)\right\}\,\right]
\end{equation}
and
\begin{equation}
\PZ=\,{1\over {2\,i}}\,(I_n-W)\PE.
\end{equation}

\vskip 0.3cm We need the following lemma for the proof of Theorem
1.2. H. Maass [3] observed the following useful fact. \vskip 0.3cm
\noindent{\bf Lemma 4.1.} (a) Let $A$ be an $m\times n$ matrix and
$B$ an $n\times l$ matrix. Assume that the entries of $A$ commute
with the entries of $B$. Then ${}^t(AB)=\,{ }^tB\,\,{ }^tA.$
\vskip 0.1cm \noindent (b) Let $A,\,B$ and $C$ be a $k\times l$,
an $n\times m$ and an $m\times l$ matrix respectively. Assume that
the entries of $A$ commute with the entries of $B$. Then

\begin{eqnarray*}
{ }^t(A\,\,{ }^t(BC))=\,B\,\,{ }^t(A\,\,^tC).
\end{eqnarray*}

\vskip 0.2cm \noindent {\it Proof.} The proof follows immediately
from the direct computation. \hfill$\Box$ \vskip 0.2cm \indent

From the formulas (3.2), (4.1) and Lemma 4.1, we get the following
formula

\begin{equation}
4\,\s\left(\,Y\,\,
^t\!\left(Y\POB\right)\PO\,\right)\,=(\alpha)+(\beta)+(\g)+(\delta),
\end{equation}
where
\begin{eqnarray*}
(\alpha):&=& \,\s \left(
(I_n-W\OW\,)\,{}^t\!\left( (I_n-W\OW\,)\PWB\right)\PW\right),\\
(\beta):&=& \,-\s \left(
\eta\,(I_n-W)^{-1}\,(I_n-W\OW\,)\,^t\!\left( (I_n-W\OW\,)\PWB\right)\PE\right),\\
(\gamma):&=&\,-\s \left( (I_n-W\OW\,)
(I_n-\OW\,)^{-1}\,{}^t\ot\,\,{}^t\!\left(\PEB\right)(I_n-\OW\,W)\PW\right),
\\
(\delta):&=&\,\s \left( \,\eta\,(I_n-W)^{-1}
(I_n-W\OW\,)(I_n-\OW)^{-1}\,{}^t\ot\,\,^t\!\left(\PEB\right)\,(I_n-\OW
W)\PE \,\right).
\end{eqnarray*}

\vskip 0.3cm According to the formulas (3.2) and (4.2), we get

\begin{equation}
4\,\s\left(\,Y\,\PZ\,\,^t\!\left(\PZB\right)\right)=\,\s\left(
(I_n-\OW W)\PE\,\,^t\!\left(\PEB\right)\right).
\end{equation}

\vskip 0.3cm From the formulas (3.2),\,(3.3) and (4.2),we get

\begin{equation}
4\,\s\left(\,VY^{-1}\,{}^tV\,\,^t\!\left(Y\PZB\right)\PZ\right)=(\epsilon)+(\zeta)+(\eta)+(\theta),
\end{equation}
where
\begin{eqnarray*}
(\epsilon):&=& \,\s \left( \eta
(I_n-\OW\,W)^{-1}\,\,{}^t\ot\,\,^t\!\left(\PEB\right)(I_n-\OW
W)\PE\right),\\
(\zeta):&=& \,\s \left( \ot\,
(I_n-W\OW\,)^{-1}\,{}^t\eta\,\,^t\!\left(\PEB\right)(I_n-\OW
W)\PE\right),\\
(\eta):&=&\,\s \left( \eta\,
(I_n-W)^{-1}(I_n-\OW\,)(I_n-W\OW\,)^{-1}\,{}^t\eta\,\,^t\!\left(\PEB\right)(I_n-\OW
W)\PE\right),
\\
(\theta):&=&\,\s \left( \,\ot\, (I_n-\OW\,)^{-1}(I_n-W)(I_n-\OW W
\,)^{-1}\,{}^t\ot\,\,^t\!\left(\PEB\right)(I_n-\OW W)\PE\right).
\end{eqnarray*}

\vskip 0.3cm Using the formulas (3.2),\,(3.3),\,(4.1),\,(4.2) and
Lemma 4.1, we get

\begin{equation}
4\,\s\left(\,V\,^t\!\left(Y\POB\right)\PZ\right)=(\iota)+(\ka)+(\la)+(\mu),
\end{equation}
where
\begin{eqnarray*}
(\iota):&=& \,\s \left( \ot \,\, ^t\!\!\left(
(I_n-W\OW\,)\PWB\right)\PE\right),\\
(\ka):&=& \,\s \left( \eta\,(I_n-W)^{-1}(I_n-\OW\,)\,^t\!\left(
(I_n-W\OW\,)\PWB\right)\PE\right),\\
(\la):&=&\,-\s \left( \eta\,
(I_n-W)^{-1}\,{}^t\ot\,\,\,^t\!\left(\PEB\right)(I_n-\OW
W)\PE\right),
\\
(\mu):&=&\,-\s \left( \,\ot\,
(I_n-\OW\,)^{-1}\,{}^t\ot\,\,\,^t\!\left(\PEB\right)(I_n-\OW
W)\PE\right).
\end{eqnarray*}

\vskip 0.3cm Using the formulas (3.2),\,(3.3),\,(4.1),\,(4.2) and
Lemma 4.1, we get

\begin{equation}
4\,\s\left(\,{}^tV\,^t\!\left(Y\PZB\right)\PO\right)=(\nu)+(\xi)+(o)+(\pi),
\end{equation}
where
\begin{eqnarray*}
(\nu):&=& \,\s \left(\,{}^t\eta \,\, ^t\!\left( \PEB
\right)(I_n-\OW W)\PW\right),\\
(\xi):&=& \,\s \left( (I_n-W)(I_n-\OW\,)^{-1}\,{}^t{\overline
\eta} \,\, ^t\!\left( \PEB
\right)(I_n-\OW W)\PW\right), \\
(o):&=&\,-\s \left( \eta\,
(I_n-W)^{-1}\,{}^t\eta\,\,\,^t\!\left(\PEB\right)(I_n-\OW
W)\PE\right),
\\
(\pi):&=&\,-\s \left( \,\eta\,
(I_n-\OW\,)^{-1}\,{}^t\ot\,\,\,^t\!\left(\PEB\right)(I_n-\OW
W)\PE\right).
\end{eqnarray*}

\vskip 0.3cm Adding the formulas $(\g),\,(\nu)$ and $(\xi)$, we
get

\begin{equation}
(\g)+(\nu)+(\xi)=\,\s\left(\,{}^t(\eta-\ot\,W)\,\,^t\!\left(\PEB\right)(I_n-\OW
W)\PW\right)
\end{equation}
because
\begin{eqnarray*}
& &\,-(I_n-W\OW\,)(I_n-\OW\,)^{-1}+(I_n-W)(I_n-\OW\,)^{-1}\\
&=&\,-W(I_n-\OW\,)(I_n-\OW\,)^{-1}=-W.
\end{eqnarray*}

\vskip 0.3cm Adding the formulas $(\beta),\,(\iota)$ and $(\ka)$,
we get
\begin{equation}
(\beta)+(\iota)+(\ka)=\,\s\left((\ot-\eta\,\OW\,)\,\,^t\!\left((I_n-W\OW\,)\PWB\right)\PE\right)
\end{equation}
because
\begin{eqnarray*}
& &\,-(I_n-W)^{-1}(I_n-W\OW\,)+(I_n-W)^{-1}(I_n-\OW\,)\\
&=&\,-(I_n-W)^{-1}(I_n-W)\OW=-\OW.
\end{eqnarray*}

\vskip 0.3cm If we add the formulas $(\eta)$ and $(o)$, we get
\begin{equation}
(\eta)+(o)=\,-\s\left(\eta\,\OW\,(I_n-W\OW\,
)^{-1}\,{}^t\eta\,\,^t\!\left(\PEB\right)(I_n-\OW W) \PE\right)
\end{equation}
because
\begin{eqnarray*}
& &\,(I_n-W)^{-1}(I_n-\OW\,)(I_n-W\OW \,)^{-1}-(I_n-W)^{-1}\\
&=&\,(I_n-W)^{-1}\Big\{
I_n-\OW-(I_n-W\OW\,)\Big\}(I_n-W\OW\,)^{-1}\\
&=&\,(I_n-W)^{-1}(I_n-W)(-\OW\,)(I_n-W\OW\,)^{-1}\\
&=&\,-\OW\,(I_n-W\OW\,)^{-1}.
\end{eqnarray*}

\vskip 0.3cm If we add the formulas $(\theta)$ and $(\mu)$, we get
\begin{equation}
(\theta)+(\mu)=\,-\s\left(\ot\,W\,(I_n-\OW
W)^{-1}\,{}^t\ot\,\,^t\!\left(\PEB\right)(I_n-\OW W) \PE\right)
\end{equation}
because
\begin{eqnarray*}
& &\,(I_n-\OW\,)^{-1}(I_n-W)(I_n-\OW W)^{-1}-(I_n-\OW\,)^{-1}\\
&=&\,(I_n-\OW\,)^{-1}\Big\{
I_n-W-(I_n-\OW W)\Big\}(I_n-\OW W)^{-1}\\
&=&\,(I_n-\OW\,)^{-1}(I_n-\OW\,)(-W)(I_n-\OW W)^{-1}\\
&=&\,-W(I_n-\OW W)^{-1}.
\end{eqnarray*}

If we add the formulas $(\delta),\,(\epsilon),\,(\la)$ and
$(\pi)$, we get

\begin{equation}
(\delta)+(\epsilon)+(\la)+(\pi)=\,\s\left(\eta\,\OW\,W\,(I_n-\OW W
)^{-1}\,{}^t\ot\,\,^t\!\left(\PEB\right)(I_n-\OW W) \PE\right)
\end{equation}
because

\begin{eqnarray*}
& &(I_n-W)^{-1}(I_n-W\OW\,)(I_n-\OW\,)^{-1}+(I_n-\OW W)^{-1}\\
& &\quad -(I_n-W)^{-1}-(I_n-\OW\,)^{-1}\\
&=&(I_n-W)^{-1}\Big\{ (I_n-W\OW\,)-(I_n-\OW\,)\Big\}(I_n-\OW\,)^{-1}\\
& &\quad +(I_n-\OW W)^{-1}-(I_n-\OW\,)^{-1}\\
&=&\,\OW\,(I_n-\OW\,)^{-1}+(I_n-\OW W)^{-1}-(I_n-\OW\,)^{-1}\\
&=&-(I_n-\OW\,)(I_n-\OW\,)^{-1}+(I_n-\OW W)^{-1}\\
&=&-I_n+(I_n-\OW W)^{-1}\\
&=&\Big\{-(I_n-\OW W)+I_n\Big\}(I_n-\OW W)^{-1}\\
&=&\OW W(I_n-\OW W)^{-1}.\\
\end{eqnarray*}

From the formulas (4.3) and (4.5)-(4.12), we obtain
\begin{eqnarray*}
&&\ \ \ \s\left(\,Y\,\,
^t\!\left(Y\POB\right)\PO\,\right)\, +\,\s\left(\,VY^{-1}\,^tV\,\,^t\!\left(Y\PZB\right)\,\PZ\,\right)\\
& &\ \ \
+\,\s\left(V\,\,^t\!\left(Y\POB\right)\PZ\,\right)+\,\s\left(\,^tV\,\,^t\!\left(Y\PZB\right)\PO\,\right)\\
&&=(\alpha)+(\beta)+(\gamma)+(\delta)+\cdots+
(\nu)+(\xi)+(o)+(\pi)\\
&&=\s\left( (I_n-W\OW)\,{}^t\!\left(
(I_n-W\OW)\PWB\right)\PW\right)\,\\
& & \quad  +\,\s \left(\,{}^t(\eta-\ot\,W)\,{}^t\!\left(
\PEB\right)
(I_n-\OW W)\PW  \right)\,\\
& &\quad  +\,\s\left( (\ot-\eta\,\OW)\,{}^t\!\left(
(I_n-W\OW)\PWB\right)\PE\right)\\
& &\quad -\,\s\left( \eta \OW
(I_n-W\OW)^{-1}\,{}^t\eta\,\,{}^t\!\left(\PEB\right)(I_n-\OW
W)\PE\right)\\
& & \quad -\,\s\left( \ot W (I_n-\OW W)^{-1}
\,{}^t\ot\,\,{}^t\!\left(\PEB\right)(I_n-\OW
W)\PE\right)\\
& &\quad  +\,\s\left( \ot (I_n-W\OW)^{-1}{}^t\eta\,{}^t\!\left(
\PEB\right)
(I_n-\OW W)\PE \right)\\
& &\quad  +\,\s\left( \eta\,\OW W (I_n-\OW
W)^{-1}\,{}^t\ot\,\,{}^t\!\left( \PEB\right) (I_n-\OW W)\PE
\right).
\end{eqnarray*}

Consequently the complete proof follows from the formula (4.4),
the above formula, Theorem B and the fact that the action (1.2) of
$G^J$ on $\Hnm$ is compatible with the action (2.6) of $G_*^J$ on
$\Dnm$ through the partial Cayley transform.\hfill $\square$

\vskip 0.3cm \noindent {\textbf{Remark 4.1.}} We proved in [11]
that the following two differential operators $D$ and $L:=\frac
14\,\Delta_{n,m;1,1}-D$ on $\Hnm$ defined by
\begin{equation*}
D=\s\left(\, Y\,\PZ\,\,{}^t\!\left( \PZB\right)\,\right)
\end{equation*}
\noindent and
\begin{eqnarray*}
L&=& \s\left(\,Y\,\,
^t\!\left(Y\POB\right)\PO\,\right)\, +\,\s\left(\,VY^{-1}\,^tV\,\,^t\!\left(Y\PZB\right)\,\PZ\,\right)\\
& &\ \
+\,\s\left(V\,\,^t\!\left(Y\POB\right)\PZ\,\right)+\,\s\left(\,^tV\,\,^t\!\left(Y\PZB\right)\PO\,\right)
\end{eqnarray*}

\noindent are invariant under the action (1.2) of $G^J$. By the
formula (4.4) and the proof of Theorem 1.2, we see that the
following differential operators ${\tilde D}$ and ${\tilde
L}:={\tilde\Delta}_{n,m;1,1}$ on $\Dnm$ defined by
\begin{equation*}
{\tilde D}=\,\s\left( (I_n-\OW W)\PE\,^t\!\left(\PEB\right)\right)
\end{equation*}
\noindent and
\begin{eqnarray*}
{\tilde L}&=& \s\left( (I_n-W\OW)\,{}^t\!\left(
(I_n-W\OW)\PWB\right)\PW\right)\,\\
& &   +\,\s \left(\,{}^t(\eta-\ot\,W)\,{}^t\!\left( \PEB\right)
(I_n-\OW W)\PW  \right)\,\\
& &  +\,\s\left( (\ot-\eta\,\OW)\,{}^t\!\left(
(I_n-W\OW)\PWB\right)\PE\right)\\
& & -\,\s\left( \eta \OW
(I_n-W\OW)^{-1}\,{}^t\eta\,\,{}^t\!\left(\PEB\right)(I_n-\OW
W)\PE\right)\\
& & -\,\s\left( \ot W (I_n-\OW W)^{-1}
\,{}^t\ot\,\,{}^t\!\left(\PEB\right)(I_n-\OW
W)\PE\right)\\
& & +\,\s\left( \ot (I_n-W\OW)^{-1}{}^t\eta\,{}^t\!\left(
\PEB\right)
(I_n-\OW W)\PE \right)\\
& &  +\,\s\left( \eta\,\OW W (I_n-\OW
W)^{-1}\,{}^t\ot\,\,{}^t\!\left( \PEB\right) (I_n-\OW W)\PE
\right)
\end{eqnarray*}

\noindent are invariant under the action (2.6) of $G_*^J.$ Indeed
it is very complicated and difficult at this moment to express the
generators of the algebra of all $G^J_{*}$-invariant differential
operators on $\Dnm$ explicitly. We propose an open problem to find
other explicit $G^J_{*}$-invariant differential operators on
$\Dnm$.

\vskip 1.05cm
\begin{section}{{\bf Remark on Harmonic Analysis on Siegel-Jacobi Disk}}
\setcounter{equation}{0}
\end{section}
\vskip 0.1cm It might be interesting to develop the theory of
harmonic analysis on the Siegel-Jacobi disk $\Dnm$. The theory of
harmonic analysis on the generalized unit disk $\BD_n$ can be done
explicitly by the work of Harish-Chandra because $\BD_n$ is a
symmetric space. However the Siegel-Jacobi disk $\Dnm$ is not a
symmetric space. The work for developing the theory of harmonic
analysis on $\Dnm$ explicitly is complicated and difficult at this
moment. We observe that this work on $\Dnm$ generalizes the work
on the generalized unit disk $\BD_n$.

\vskip 2mm More precisely, if we put $G_*=SU(n,n)\cap Sp(n,\BC)$,
then the Jacobi group
\begin{equation*}
G_*^J=\bigg\{ \left( \begin{pmatrix} P & Q \\ {\overline Q} &
{\overline P} \end{pmatrix},(\xi,{\overline \xi};\,i\kappa)
\right)\,\Big|\ \begin{pmatrix} P & Q \\ {\overline Q} &
{\overline P} \end{pmatrix}\in G_*,\ \xi\in\BC^{(m,n)},\ \kappa\in
\BR^{(m,m)}\,\bigg\}
\end{equation*}

\noindent acts on the Siegel-Jacobi disk $\Dnm$ transitively via
the transformation behavior (1.9). It is easily seen that the
stabilizer $K_*^J$ of the action (1.9) at the base point $(0,0)$
is given by
\begin{equation*}
K_*^J=\bigg\{ \left( \begin{pmatrix} P & 0 \\ 0 & {\overline P}
\end{pmatrix},(0,0;\,i\kappa) \right)\,\Big|\
P\in U(n),\ \kappa\in \BR^{(m,m)}\,\bigg\}.
\end{equation*}

\noindent Therefore $G_*^J/K_*^J$ is biholomorphic to $\Dnm$ via
the correspondence
\begin{equation*}
gK_*^J\mapsto g\cdot (0,0),\quad g\in G_*^J.
\end{equation*}

\noindent We observe that the Siegel-Jacobi disk $\Dnm$ is not a
reductive symmetric space.

\vskip 2mm We let
\begin{equation*}
\G_{n,m}:=Sp(n,{\mathbb Z})\ltimes H_{\mathbb Z}^{(n,m)},
\end{equation*}

\noindent where $Sp(n,\mathbb Z)$ is the Siegel modular group of
degree $n$ and
\begin{equation*}
H_{\mathbb Z}^{(n,m)}=\Big\{ (\lambda,\mu;\kappa)\in
H_{\BR}^{(n,m)}\,|\ \lambda,\mu,\kappa \ \textrm{are
integral}\,\Big\}.
\end{equation*}

\noindent We set
\begin{equation*}
\G_{n,m}^*:=T_*^{-1}\G_{n,m}T_*,
\end{equation*}

\noindent where $T_*$ was already defined in Section 2. Clearly
the arithmetic subgroup $\G_{n,m}^*$ acts on $\Dnm$ properly
continuously. We can describe a fundamental domain ${\mathcal
F}_{n,m}^*$ for $\G_{n,m}^*\backslash \Dnm$ explicitly using the
partial Cayley transform and a fundamental domain ${\mathcal
F}_{n,m}$ for $\G_{n,m}\backslash \Hnm$ which is described
explicitly in [13]. The $G_*^J$-invariant metric $d{\tilde
s}_{n,m;A,B}$ on $\Dnm$ induces a metric on ${\mathcal F}_{n,m}^*$
naturally. It may be intersting to investigate the spectral theory
of the Laplacian ${\tilde \Delta}_{n,m;A,B}$ on a fundamental
domain ${\mathcal F}_{n,m}^*$. But this work is very complicated
and difficult at this moment.

\vskip 2mm For instance, we consider the case $n=m=1$ and $A=B=1$.
In this case
\begin{equation*}
G_*^J=\bigg\{ \left( \begin{pmatrix} p & q \\ {\overline q} &
{\overline p} \end{pmatrix},(\xi,{\overline \xi};\,i\kappa)
\right)\,\Big|\ p,q,\xi\in \BC,\ |p|^2-|q|^2=1,\ \kappa\in
\BR\,\bigg\}
\end{equation*}

\noindent and
\begin{equation*}
K_*^J=\bigg\{ \left( \begin{pmatrix} p & 0 \\ 0 & {\overline p}
\end{pmatrix},(0,0;\,i\kappa) \right)\,\Big|\
p\in\BC,\ |p|=1,\ \kappa\in \BR\,\bigg\}.
\end{equation*}

\noindent $d{\tilde s}_{1,1;1,1}$ is a $G_*^J$-invariant
Riemannian metric on $\BD_{1,1}=\BD_1\times \BC$\,(cf. Theorem
1.1) and ${\tilde \Delta}_{1,1;1,1}$ is its Laplacian. It is well
known that the theory of harmonic analysis on the unit disk
$\BD_1$ has been well developed explicitly (cf.\,[3], pp.\,29-72).
I think that so far nobody has not investigated the theory of
harmonic analysis on $\BD_{1,1}$ explicitly. For example,
inversion formula, Plancherel formula, Paley-Wiener theorem on
$\BD_{1,1}$ have not been described explicitly until now. It seems
that it is interesting to develop the theory of harmonic analysis
on the Siegel-Jacobi disk $\BD_{1,1}$ explicitly.

\end{document}